\documentclass[10pt]{article}
\usepackage{cite}
\usepackage{mathrsfs}
\usepackage{amsfonts}
\usepackage{amsmath}
\usepackage{amsfonts,amssymb}
\usepackage{dsfont}
\usepackage{curves}
\usepackage{mathrsfs}
\usepackage{pifont}
\usepackage{amssymb}
\allowdisplaybreaks

\numberwithin{equation}{section}

\date{}

\def\BigRoman{\uppercase\expandafter{\romannumeral\number\count 255 }}
\def\Romannumeral{\afterassignment\BigRoman\count255=}

\begin{document}
\title{Some sufficient conditions for graphs to have component factors
}
\author{\small  Sizhong Zhou\footnote{Corresponding author. E-mail address: zsz\_cumt@163.com (S. Zhou)}\\
\small  School of Science, Jiangsu University of Science and Technology,\\
\small  Zhenjiang, Jiangsu 212100, China\\
}

\maketitle
\begin{abstract}
\noindent Let $G$ denote a graph and $k\geq2$ be an integer. A $\{K_{1,1},K_{1,2},\ldots,K_{1,k},\mathcal{T}(2k+1)\}$-factor of $G$ is a spanning
subgraph, whose every connected component is isomorphic to an element of $\{K_{1,1},K_{1,2},\ldots,K_{1,k},\mathcal{T}(2k+1)\}$, where
$\mathcal{T}(2k+1)$ is one special family of tree. In this paper, we put forward some sufficient conditions for the existence of
$\{K_{1,1},K_{1,2},\ldots,K_{1,k},\mathcal{T}(2k+1)\}$-factors in graphs. Furthermore, we construct some extremal graphs to show that the main
results in this paper are best possible.
\\
\begin{flushleft}
{\em Keywords:}  $\{K_{1,1},K_{1,2},\ldots,K_{1,k},\mathcal{T}(2k+1)\}$-factor; Laplacian eigenvalue; size; degree; independence number.

(2020) Mathematics Subject Classification: 05C70, 05C50
\end{flushleft}
\end{abstract}

\section{Introduction}

In this paper, we deal only with finite and undirected graphs without loops and multiple edges. Let $G$ denote a graph with vertex set $V(G)$
and edge set $E(G)$. The order of $G$ is the number $n=|V(G)|$ of its vertices and its size is the number $m=|E(G)|$ of its edges. A graph of
order 1 is called trivial. We use $I(G)$ to denote the set of isolated vertices in $G$ and write $i(G)=|I(G)|$. The degree of a vertex $x$ in
$G$, denote by $d_G(x)$, is defined as the number of edges which are incident to $x$. Let $\omega(G)$, $\delta(G)$ and $\alpha(G)$ denote the
number of connected components, the minimum degree and the independence number of $G$, respectively. For any $S\subseteq V(G)$, we denote by
$G[S]$ the subgraph of $G$ induced by $S$, and by $G-S$ the graph formed from $G$ by deleting the vertices in $S$ and their incident edges.
Let $c$ be a real number. As usual, let $C_n$, $P_n$, $K_{1,n-1}$ and $K_n$ denote the cycle, the path, the star and the complete graph of
order $n$, respectively. Recall that $\lceil c\rceil$ is the least integer satisfying $\lceil c\rceil\geq c$ and $\lfloor c\rfloor$ is the
greatest integer satisfying $\lfloor c\rfloor\leq c$.

We use $Leaf(T)$ to denote the set of leaves in a tree $T$. An edge of a tree $T$ incident with a leaf is called a pendant edge. The number of
pendant edges in a tree $T$ is equal to that of leaves in $T$ under the case that the order of $T$ is at least 3. We are to define a special
class of trees $\mathcal{T}(2k+1)$, where $k\geq2$ is an integer. Let $R$ be a tree that satisfies the following conditions: for any
$x\in V(R)-Leaf(R)$,

(a) $d_{R-Leaf(R)}(x)\in\{1,3,\ldots,2k+1\}$\\
and

(b) 2(the number of leaves adjacent to $x$ in $R$)$+d_{R-Leaf(R)}(x)\leq2k+1$.\\
For such a tree $R$, we derive a new tree $T_R$ as follows:

(c) insert a new vertex of degree 2 into each edge of $R-Leaf(R)$\\
and

(d) for every vertex $x$ of $R-Leaf(R)$ with $d_{R-Leaf(R)}(x)=2r+1<2k+1$, add $k-r-$(the number of leaves adjacent to $x$ in $R$)
pendant edges to $x$.\\
Then the set of such trees $T_R$ for all trees $R$ satisfying conditions (a) and (b) is denoted by $\mathcal{T}(2k+1)$.

Let $G_1$ and $G_2$ be two vertex-disjoint graphs. Let $G_1\cup G_2$ denote the union of $G_1$ and $G_2$. The join $G_1\vee G_2$ is obtained
from $G_1\cup G_2$ by joining each vertex of $G_1$ with each vertex of $G_2$ by an edge.

Given a graph $G$ of order $n$, the adjacency matrix $A(G)=(a_{ij})_{n\times n}$ of $G$ is a 0--1 matrix in which the entry $a_{ij}=1$ if and
only if $v_i$ and $v_j$ are adjacent. Let $D(G)$ be the diagonal matrix of vertex degrees of $G$. Let $L(G)=D(G)-A(G)$ be the Laplacian matrix
of $G$. Let $\mu_1(G)\geq\mu_2(G)\geq\cdots\geq\mu_n(G)=0$ be eigenvalues of $L(G)$.

Let $\mathcal{H}$ denote a set of connected graphs. Then a spanning subgraph $H$ of $G$ is called an $\mathcal{H}$-factor if each connected
component of $H$ is isomorphic to an element of $\mathcal{H}$. An $\mathcal{H}$-factor is also referred as a component factor. Let $k\geq2$
be an integer. Write $P_{\geq k}=\{P_i:i\geq k\}$. If $\mathcal{H}=\{P_i:i\geq k\}$, then an $\mathcal{H}$-factor is called a
$P_{\geq k}$-factor. If $\mathcal{H}=\{K_{1,1},K_{1,2},\ldots,K_{1,k},\mathcal{T}(2k+1)\}$, then an $\mathcal{H}$-factor is called a
$\{K_{1,1},K_{1,2},\ldots,K_{1,k},\mathcal{T}(2k+1)\}$-factor. In particular, a perfect matching is also a $\{K_{1,1}\}$-factor.

Li and Miao \cite{LM}, Zhou, Zhang and Sun \cite{ZZS}, Zhou, Sun and Liu \cite{ZSL1} established some spectral radius conditions for graphs to
contain $P_{\geq2}$-factors. Many researchers \cite{GCW,WZi,Wp,LP,ZSL2,Zd,Zp,Zs,ZSB} proved some results related to binding number, neighborhood,
degree condition etc., for a graph to contain a $P_{\geq3}$-factor. O \cite{Os} established a connection between spectral radius and
$\{K_{1,1}\}$-factors in graphs. Tutte \cite{T} provided a characterization for a graph with a $\{K_2,C_n:n\geq3\}$-factor. Klopp and Steffen
\cite{KSf} studied the existence of a $\{K_{1,1},K_{1,2},C_m:m\geq3\}$-factor in a graph. Amahashi and Kano \cite{AK} showed a necessary and
sufficient condition for a graph to have a $\{K_{1,j}:1\leq j\leq k\}$-factor, where $k$ is an integer with $k\geq2$. Kano and Saito \cite{KS}
presented a sufficient condition for a graph to have a $\{K_{1,j}:k\leq j\leq2k\}$-factor, where $k\geq2$ is an integer. Kano, Lu and Yu
\cite{KLYc} showed a relationship between the number of isolated vertices and $\{K_{1,2},K_{1,3},K_5\}$-factors in graphs. Many efforts have
been devoted to finding sufficient conditions for graphs to contain spanning subgraphs by utilizing various graphic parameters such as spectral
radius \cite{ZSL3,ZZB,ZZL,Wc}, independence number \cite{KL,Wa}, degree condition\cite{Lv,ZPX1}, isolated toughness \cite{ZPX2}, eigenvalues
\cite{KOPR,Oe} and others \cite{Za1,GWW}.

Kano, Lu and Yu \cite{KLYf} provided a necessary and sufficient condition for a graph to have a $\{K_{1,1},K_{1,2},\\
\cdots,K_{1,k},\mathcal{T}(2k+1)\}$-factor.

\medskip

\noindent{\textbf{Theorem 1.1}} (Kano, Lu and Yu \cite{KLYf}). Let $k\geq2$ be an integer. Then a graph $G$ contains a
$\{K_{1,1},K_{1,2},\ldots,K_{1,k},\mathcal{T}(2k+1)\}$-factor if and only if
$$
i(G-S)\leq\Big(k+\frac{1}{2}\Big)|S|
$$
for any $S\subseteq V(G)$.

\medskip

Motivated by \cite{KLYf}, it is natural and interesting to consider the existence of $\{K_{1,1},K_{1,2},\ldots,K_{1,k},\mathcal{T}(2k+1)\}$-factor
in graph. Here, we focus on the sufficient conditions including Laplacian eigenvalue condition, size condition, degree condition or independence
number condition. Our main results are shown in the following.

\medskip

\noindent{\textbf{Theorem 1.2.}} Let $k\geq2$ be an integer, and let $G$ be a graph of order $n$. If
$$
\mu_1\leq\Big(k+\frac{1}{2}\Big)\mu_{n-1},
$$
then $G$ has a $\{K_{1,1},K_{1,2},\ldots,K_{1,k},\mathcal{T}(2k+1)\}$-factor.

\medskip

\noindent{\textbf{Theorem 1.3.}} Let $k$ and $t$ be two integers with $k\geq2$ and $1\leq t\leq k-1$, and let $G$ be a $t$-connected graph of order
$n$ with $n\geq\frac{(2k^{2}+5k+1)t+k^{2}+5k+2}{2k}$. If $G$ satisfies
$$
|E(G)|>\binom{n-\left\lfloor\left(k+\frac{1}{2}\right)t\right\rfloor-1}{2}+t\left(\left\lfloor\left(k+\frac{1}{2}\right)t\right\rfloor+1\right),
$$
then $G$ has a $\{K_{1,1},K_{1,2},\ldots,K_{1,k},\mathcal{T}(2k+1)\}$-factor.
\medskip

\noindent{\textbf{Theorem 1.4.}} Let $k\geq2$ be an integer, and let $G$ be a graph of order $n$ with minimum degree $\delta\geq1$. If $G$ satisfies
$$
\max\{d_G(x_1),d_G(x_2),\ldots,d_G(x_{\lfloor(k+\frac{1}{2})\delta\rfloor+1})\}\geq\frac{2n}{2k+3}
$$
for any independent subset $\{x_1,x_2,\ldots,x_{\lfloor(k+\frac{1}{2})\delta\rfloor+1}\}$ of $G$, then $G$ has a $\{K_{1,1},K_{1,2},\ldots,K_{1,k},\mathcal{T}(2k+1)\}$-factor.

\medskip

\noindent{\textbf{Theorem 1.5.}} Let $k\geq2$ be an integer, and let $G$ be a graph. If $G$ satisfies
$$
\delta(G)\geq\frac{2\alpha(G)}{2k+1},
$$
then $G$ has a $\{K_{1,1},K_{1,2},\ldots,K_{1,k},\mathcal{T}(2k+1)\}$-factor.

\medskip

\section{The proof of Theorem 1.2}

In order to prove Theorem 1.2, we first pose the following lemma which is obtained by Gu and Liu \cite{GL}.

\medskip

\noindent{\textbf{Lemma 2.1}} (Gu and Liu \cite{GL}). Let $G$ be a graph with $n$ vertices and at least one edge. Suppose that $S\subset V(G)$ such
that $G-S$ is disconnected. Let $X$ and $Y$ be disjoint vertex subsets of $G-S$ such that $X\cup Y=V(G)-S$ with $|X|\leq|Y|$. Then
$$
|X|\leq\frac{\mu_1-\mu_{n-1}}{2\mu_1}n,
$$
and
$$
|S|\geq\frac{2\mu_{n-1}}{\mu_1-\mu_{n-1}}|X|,
$$
with each equality holding only when $|X|=|Y|$.

\medskip

In what follows, we verify Theorem 1.2.

\medskip

\noindent{\it Proof of Theorem 1.2.} Suppose, to the contrary, that $G$ has no $\{K_{1,1},K_{1,2},\ldots,K_{1,k},\mathcal{T}(2k+1)\}$-factor. Then by
Theorem 1.1, we obtain
\begin{align}\label{eq:2.1}
i(G-S)>\left(k+\frac{1}{2}\right)|S|
\end{align}
for some subset $S$ of $V(G)$. Let $G_1,G_2,\ldots,G_{\omega}$ be components of $G-S$ with $|V(G_1)|\leq|V(G_2)|\leq\cdots\leq|V(G_{\omega})|$, where
$\omega=\omega(G-S)$. Obviously, $\omega=\omega(G-S)\geq i(G-S)=i$. If $\omega>i$, let $X=\cup_{1\leq j\leq\lceil\frac{i+1}{2}\rceil}V(G_j)$ and
$Y=V(G)-S-X$. Then we have $\frac{i+1}{2}\leq|X|\leq|Y|$. By virtue of Lemma 2.1 and $\mu_1\leq(k+\frac{1}{2})\mu_{n-1}$, we deduce
\begin{align*}
|S|\geq&\frac{2\mu_{n-1}}{\mu_1-\mu_{n-1}}|X|\geq\frac{2\mu_{n-1}}{(k+\frac{1}{2})\mu_{n-1}-\mu_{n-1}}|X|\\
=&\frac{2}{k-\frac{1}{2}}|X|\geq\frac{2}{k-\frac{1}{2}}\cdot\frac{i+1}{2}=\frac{2(i+1)}{2k-1},
\end{align*}
which leads to $i(G-S)=i\leq(k-\frac{1}{2})|S|-1=(k+\frac{1}{2})|S|-|S|-1\leq(k+\frac{1}{2})|S|-1$, which is a contradiction to \eqref{eq:2.1}. If
$\omega=i$ and $i$ is odd, we let $X=\cup_{1\leq j\leq\frac{i-1}{2}}V(G_j)$ and $Y=V(G)-S-X$. Then we conclude $|X|=\frac{i-1}{2}<\frac{i+1}{2}=|Y|$.
Together with Lemma 2.1 and $\mu_1\leq(k+\frac{1}{2})\mu_{n-1}$, we obtain
\begin{align*}
|S|>&\frac{2\mu_{n-1}}{\mu_1-\mu_{n-1}}|X|\geq\frac{2\mu_{n-1}}{(k+\frac{1}{2})\mu_{n-1}-\mu_{n-1}}|X|\\
=&\frac{2}{k-\frac{1}{2}}\cdot\frac{i-1}{2}=\frac{i-1}{k-\frac{1}{2}},
\end{align*}
which yields
\begin{align}\label{eq:2.2}
i(G-S)=i<\left(k-\frac{1}{2}\right)|S|+1=\left(k+\frac{1}{2}\right)|S|-|S|+1.
\end{align}
If $S=\emptyset$, then it follows from \eqref{eq:2.1} and \eqref{eq:2.2} that $0<i(G)<1$, which is impossible since $i(G)$ is nonnegative integer.
If $S\neq\emptyset$, then by \eqref{eq:2.2} we have $i(G-S)=i<\left(k+\frac{1}{2}\right)|S|-|S|+1\leq\left(k+\frac{1}{2}\right)|S|$, which is a
contradiction to \eqref{eq:2.1}. If $\omega=i$ and $i$ is even, we let $X=\cup_{1\leq j\leq\frac{i}{2}}V(G_j)$ and $Y=V(G)-S-X$. Then we obtain
$|X|=|Y|=\frac{i}{2}$. According to Lemma 2.1 and $\mu_1\leq(k+\frac{1}{2})\mu_{n-1}$, we get
\begin{align*}
|S|\geq&\frac{2\mu_{n-1}}{\mu_1-\mu_{n-1}}|X|\geq\frac{2\mu_{n-1}}{(k+\frac{1}{2})\mu_{n-1}-\mu_{n-1}}|X|\\
=&\frac{2}{k-\frac{1}{2}}|X|=\frac{2i}{2k-1},
\end{align*}
which implies $i(G-S)=i\leq(k-\frac{1}{2})|S|=(k+\frac{1}{2})|S|-|S|\leq(k+\frac{1}{2})|S|$, which contradicts \eqref{eq:2.1}. This completes the
proof of Theorem 1.2. \hfill $\Box$

\section{The proof of Theorem 1.3}

\noindent{\it Proof of Theorem 1.3.} Suppose, to the contrary, that $G$ contains no $\{K_{1,1},K_{1,2},\ldots,K_{1,k},\mathcal{T}(2k+1)\}$-factor.
Then by Theorem 1.1, we possess
$$
i(G-S)>\left(k+\frac{1}{2}\right)|S|
$$
for some subset $S$ of $V(G)$. According to the integrity of $i(G-S)$, we obtain
\begin{align}\label{eq:3.1}
i(G-S)\geq\left\lfloor\left(k+\frac{1}{2}\right)|S|\right\rfloor+1
\end{align}
for some subset $S$ of $V(G)$. Obviously, $|S|\geq t$. Otherwise, $G-S$ is connected and so $i(G-S)=0$. Together with \eqref{eq:3.1}, we get
$0=i(G-S)\geq\left\lfloor\left(k+\frac{1}{2}\right)|S|\right\rfloor+1\geq1$, a contradiction. Let $|S|=s\geq t$. Then $G$ is a spanning subgraph
of $G_1=K_s\vee(K_{n_1}\cup\left(\left\lfloor\left(k+\frac{1}{2}\right)s\right\rfloor+1\right)K_1)$, where
$n_1=n-\left\lfloor\left(k+\frac{3}{2}\right)s\right\rfloor-1\geq0$ is an integer. Thus, we get
\begin{align}\label{eq:3.2}
|E(G)|\leq |E(G_1)|.
\end{align}
Notice that $n\geq\left\lfloor\left(k+\frac{3}{2}\right)s\right\rfloor+1$ and
$|E(G_1)|=\binom{n-\left\lfloor\left(k+\frac{1}{2}\right)s\right\rfloor-1}{2}+s\left(\left\lfloor\left(k+\frac{1}{2}\right)s\right\rfloor+1\right)$.
If $s=t$, then $|E(G_1)|=\binom{n-\left\lfloor\left(k+\frac{1}{2}\right)t\right\rfloor-1}{2}+t\left(\left\lfloor\left(k+\frac{1}{2}\right)t\right\rfloor+1\right)$.
Combining this with \eqref{eq:3.2}, we conclude
$|E(G)|\leq |E(G_1)|=\binom{n-\left\lfloor\left(k+\frac{1}{2}\right)t\right\rfloor-1}{2}+t\left(\left\lfloor\left(k+\frac{1}{2}\right)t\right\rfloor+1\right)$,
a contradiction. In what follows, we consider $s\geq t+1$.

\noindent{\bf Case 1.} $s$ and $t$ are odd.

In this case, $s\geq t+2$, $n\geq\left\lfloor\left(k+\frac{3}{2}\right)s\right\rfloor+1=\left(k+\frac{3}{2}\right)s+\frac{1}{2}$,
$\left\lfloor\left(k+\frac{1}{2}\right)t\right\rfloor=\left(k+\frac{1}{2}\right)t-\frac{1}{2}$ and
$\left\lfloor\left(k+\frac{1}{2}\right)s\right\rfloor=\left(k+\frac{1}{2}\right)s-\frac{1}{2}$. By a direct computation, we obtain
\begin{align*}
&\binom{n-\left\lfloor\left(k+\frac{1}{2}\right)t\right\rfloor-1}{2}+t\left(\left\lfloor\left(k+\frac{1}{2}\right)t\right\rfloor+1\right)-|E(G_1)|\\
=&\binom{n-\left(k+\frac{1}{2}\right)t-\frac{1}{2}}{2}+t\left(\left(k+\frac{1}{2}\right)t+\frac{1}{2}\right)-
\binom{n-\left(k+\frac{1}{2}\right)s-\frac{1}{2}}{2}-s\left(\left(k+\frac{1}{2}\right)s+\frac{1}{2}\right)\\
=&\frac{1}{8}(s-t)((8k+4)n-(2k+1)(2k+5)(s+t)-8k-8)\\
\geq&\frac{1}{8}(s-t)\left((8k+4)\left(\left(k+\frac{3}{2}\right)s+\frac{1}{2}\right)-(2k+1)(2k+5)(s+t)-8k-8\right)\\
=&\frac{1}{8}(s-t)((2k+1)^{2}s-(2k+1)(2k+5)t-4k-6)\\
\geq&\frac{1}{8}(s-t)((2k+1)^{2}(t+2)-(2k+1)(2k+5)t-4k-6)\\
=&\frac{1}{8}(s-t)(-4(2k+1)t+2(2k+1)^{2}-4k-6)\\
\geq&\frac{1}{8}(s-t)(-4(2k+1)(k-1)+2(2k+1)^{2}-4k-6)\\
=&k(s-t)\\
>&0
\end{align*}
by $s\geq t+2$, $n\geq\left(k+\frac{3}{2}\right)s+\frac{1}{2}$ and $1\leq t\leq k-1$, which implies
\begin{align}\label{eq:3.3}
|E(G_1)|<\binom{n-\left\lfloor\left(k+\frac{1}{2}\right)t\right\rfloor-1}{2}+t\left(\left\lfloor\left(k+\frac{1}{2}\right)t\right\rfloor+1\right).
\end{align}

According to \eqref{eq:3.2} and \eqref{eq:3.3}, we infer
$$
|E(G)|\leq |E(G_1)|<\binom{n-\left\lfloor\left(k+\frac{1}{2}\right)t\right\rfloor-1}{2}+t\left(\left\lfloor\left(k+\frac{1}{2}\right)t\right\rfloor+1\right),
$$
which contradicts $|E(G)|>\binom{n-\left\lfloor\left(k+\frac{1}{2}\right)t\right\rfloor-1}{2}+t\left(\left\lfloor\left(k+\frac{1}{2}\right)t\right\rfloor+1\right)$.

\noindent{\bf Case 2.} $s$ and $t$ are even.

In this case, $s\geq t+2$, $n\geq\left\lfloor\left(k+\frac{3}{2}\right)s\right\rfloor+1=\left(k+\frac{3}{2}\right)s+1$,
$\left\lfloor\left(k+\frac{1}{2}\right)t\right\rfloor=\left(k+\frac{1}{2}\right)t$ and
$\left\lfloor\left(k+\frac{1}{2}\right)s\right\rfloor=\left(k+\frac{1}{2}\right)s$. By a direct calculation, we get
\begin{align*}
&\binom{n-\left\lfloor\left(k+\frac{1}{2}\right)t\right\rfloor-1}{2}+t\left(\left\lfloor\left(k+\frac{1}{2}\right)t\right\rfloor+1\right)-|E(G_1)|\\
=&\binom{n-\left(k+\frac{1}{2}\right)t-1}{2}+t\left(\left(k+\frac{1}{2}\right)t+1\right)-
\binom{n-\left(k+\frac{1}{2}\right)s-1}{2}-s\left(\left(k+\frac{1}{2}\right)s+1\right)\\
=&\frac{1}{8}(s-t)((8k+4)n-(2k+1)(2k+5)(s+t)-12k-14)\\
\geq&\frac{1}{8}(s-t)\left((8k+4)\left(\left(k+\frac{3}{2}\right)s+1\right)-(2k+1)(2k+5)(s+t)-12k-14\right)\\
=&\frac{1}{8}(s-t)((2k+1)^{2}s-(2k+1)(2k+5)t-4k-10)\\
\geq&\frac{1}{8}(s-t)((2k+1)^{2}(t+2)-(2k+1)(2k+5)t-4k-10)\\
=&\frac{1}{2}(s-t)(-(2k+1)t+k(2k+1)-2)\\
\geq&\frac{1}{2}(s-t)(-(2k+1)(k-1)+k(2k+1)-2)\\
=&\frac{1}{2}(s-t)(2k-1)\\
>&0
\end{align*}
by $s\geq t+2$, $n\geq\left(k+\frac{3}{2}\right)s+1$ and $1\leq t\leq k-1$, which yields
\begin{align*}
|E(G_1)|<\binom{n-\left\lfloor\left(k+\frac{1}{2}\right)t\right\rfloor-1}{2}+t\left(\left\lfloor\left(k+\frac{1}{2}\right)t\right\rfloor+1\right).
\end{align*}
Combining this with \eqref{eq:3.2}, we conclude
$$
|E(G)|\leq |E(G_1)|<\binom{n-\left\lfloor\left(k+\frac{1}{2}\right)t\right\rfloor-1}{2}+t\left(\left\lfloor\left(k+\frac{1}{2}\right)t\right\rfloor+1\right),
$$
which contradicts $|E(G)|>\binom{n-\left\lfloor\left(k+\frac{1}{2}\right)t\right\rfloor-1}{2}+t\left(\left\lfloor\left(k+\frac{1}{2}\right)t\right\rfloor+1\right)$.

\noindent{\bf Case 3.} $s$ is odd and $t$ is even.

In this case, $s\geq t+1$, $n\geq\left\lfloor\left(k+\frac{3}{2}\right)s\right\rfloor+1=\left(k+\frac{3}{2}\right)s+\frac{1}{2}$,
$\left\lfloor\left(k+\frac{1}{2}\right)t\right\rfloor=\left(k+\frac{1}{2}\right)t$ and
$\left\lfloor\left(k+\frac{1}{2}\right)s\right\rfloor=\left(k+\frac{1}{2}\right)s-\frac{1}{2}$. By a direct computation, we have
\begin{align}\label{eq:3.4}
&\binom{n-\left\lfloor\left(k+\frac{1}{2}\right)t\right\rfloor-1}{2}+t\left(\left\lfloor\left(k+\frac{1}{2}\right)t\right\rfloor+1\right)-|E(G_1)|\nonumber\\
=&\binom{n-\left(k+\frac{1}{2}\right)t-1}{2}+t\left(\left(k+\frac{1}{2}\right)t+1\right)-
\binom{n-\left(k+\frac{1}{2}\right)s-\frac{1}{2}}{2}-s\left(\left(k+\frac{1}{2}\right)s+\frac{1}{2}\right)\nonumber\\
=&\frac{1}{8}((2(2k+1)(s-t)-2)\cdot2n-(2k+1)(2k+5)s^{2}-(8k+8)s\nonumber\\
&+(2k+1)(2k+5)t^{2}+(12k+14)t+5).
\end{align}

If $s=t+1$, then it follows from \eqref{eq:3.4} and $n\geq\frac{(2k^{2}+5k+1)t+k^{2}+5k+2}{2k}$ that
\begin{align*}
&\binom{n-\left\lfloor\left(k+\frac{1}{2}\right)t\right\rfloor-1}{2}+t\left(\left\lfloor\left(k+\frac{1}{2}\right)t\right\rfloor+1\right)-|E(G_1)|\\
=&\frac{1}{8}(8kn-(2k+1)(2k+5)(t+1)^{2}-(8k+8)(t+1)+(2k+1)(2k+5)t^{2}+(12k+14)t+5)\\
=&\frac{1}{2}(2kn-(2k^{2}+5k+1)t-k^{2}-5k-2)\\
\geq&0,
\end{align*}
which leads to $|E(G_1)|\leq\binom{n-\left\lfloor\left(k+\frac{1}{2}\right)t\right\rfloor-1}{2}+t\left(\left\lfloor\left(k+\frac{1}{2}\right)t\right\rfloor+1\right)$.
Combining this with \eqref{eq:3.2}, we conclude
$|E(G)|\leq |E(G_1)|\leq\binom{n-\left\lfloor\left(k+\frac{1}{2}\right)t\right\rfloor-1}{2}+t\left(\left\lfloor\left(k+\frac{1}{2}\right)t\right\rfloor+1\right)$,
a contradiction. If $s\geq t+2$, then it follows from \eqref{eq:3.4} and $n\geq\left(k+\frac{3}{2}\right)s+\frac{1}{2}$ that
\begin{align}\label{eq:3.5}
&\binom{n-\left\lfloor\left(k+\frac{1}{2}\right)t\right\rfloor-1}{2}+t\left(\left\lfloor\left(k+\frac{1}{2}\right)t\right\rfloor+1\right)-|E(G_1)|\nonumber\\
\geq&\frac{1}{8}((2(2k+1)(s-t)-2)\cdot((2k+3)s+1)-(2k+1)(2k+5)s^{2}-(8k+8)s\nonumber\\
&+(2k+1)(2k+5)t^{2}+(12k+14)t+5)\nonumber\\
=&\frac{1}{8}((2k+1)^{2}s^{2}-(2(2k+1)(2k+3)t+8k+12)s+(2k+1)(2k+5)t^{2}+(8k+12)t+3)\nonumber\\
=&\frac{1}{8}\varphi(s),
\end{align}
where $\varphi(s)=(2k+1)^{2}s^{2}-(2(2k+1)(2k+3)t+8k+12)s+(2k+1)(2k+5)t^{2}+(8k+12)t+3$. Note that
$$
\frac{2(2k+1)(2k+3)t+8k+12}{2(2k+1)^{2}}<t+2\leq s
$$
by $1\leq t\leq k-1$. Then we have
\begin{align*}
\varphi(s)\geq&\varphi(t+2)\\
=&-8(2k+1)t+4(2k+1)^{2}-16k-21\\
\geq&-8(2k+1)(k-1)+4(2k+1)^{2}-16k-21\\
=&8k-9\\
>&0.
\end{align*}
Combining this with \eqref{eq:3.5}, we deduce
$$
\binom{n-\left\lfloor\left(k+\frac{1}{2}\right)t\right\rfloor-1}{2}+t\left(\left\lfloor\left(k+\frac{1}{2}\right)t\right\rfloor+1\right)-|E(G_1)|
\geq\frac{1}{8}\varphi(s)>0,
$$
which implies $|E(G_1)|<\binom{n-\left\lfloor\left(k+\frac{1}{2}\right)t\right\rfloor-1}{2}+t\left(\left\lfloor\left(k+\frac{1}{2}\right)t\right\rfloor+1\right)$.
Combining this with \eqref{eq:3.2}, we conclude
$|E(G)|\leq |E(G_1)|<\binom{n-\left\lfloor\left(k+\frac{1}{2}\right)t\right\rfloor-1}{2}+t\left(\left\lfloor\left(k+\frac{1}{2}\right)t\right\rfloor+1\right)$,
a contradiction.

\noindent{\bf Case 4.} $s$ is even and $t$ is odd.

In this case, $s\geq t+1$, $n\geq\left\lfloor\left(k+\frac{3}{2}\right)s\right\rfloor+1=\left(k+\frac{3}{2}\right)s+1$,
$\left\lfloor\left(k+\frac{1}{2}\right)t\right\rfloor=\left(k+\frac{1}{2}\right)t-\frac{1}{2}$ and
$\left\lfloor\left(k+\frac{1}{2}\right)s\right\rfloor=\left(k+\frac{1}{2}\right)s$. By a direct calculation, we obtain
\begin{align}\label{eq:3.6}
&\binom{n-\left\lfloor\left(k+\frac{1}{2}\right)t\right\rfloor-1}{2}+t\left(\left\lfloor\left(k+\frac{1}{2}\right)t\right\rfloor+1\right)-|E(G_1)|\nonumber\\
=&\binom{n-\left(k+\frac{1}{2}\right)t-\frac{1}{2}}{2}+t\left(\left(k+\frac{1}{2}\right)t+\frac{1}{2}\right)-
\binom{n-\left(k+\frac{1}{2}\right)s-1}{2}-s\left(\left(k+\frac{1}{2}\right)s+1\right)\nonumber\\
=&\frac{1}{8}((2(2k+1)(s-t)+2)\cdot2n-(2k+1)(2k+5)s^{2}-(12k+14)s\nonumber\\
&+(2k+1)(2k+5)t^{2}+(8k+8)t-5).
\end{align}
If $s=t+1$, then it follows from \eqref{eq:3.6} and $n\geq\frac{(2k^{2}+5k+1)t+k^{2}+5k+2}{2k}$ that
\begin{align*}
&\binom{n-\left\lfloor\left(k+\frac{1}{2}\right)t\right\rfloor-1}{2}+t\left(\left\lfloor\left(k+\frac{1}{2}\right)t\right\rfloor+1\right)-|E(G_1)|\\
=&\frac{1}{8}((8k+8)n-(2k+1)(2k+5)(t+1)^{2}-(12k+14)(t+1)+(2k+1)(2k+5)t^{2}+(8k+8)t-5)\\
=&\frac{1}{2}((2k+2)n-(2k^{2}+7k+4)t-k^{2}-6k-6)\\
\geq&\frac{1}{2}\left((2k+2)\cdot\frac{(2k^{2}+5k+1)t+k^{2}+5k+2}{2k}-(2k^{2}+7k+4)t-k^{2}-6k-6\right)\\
=&\frac{(2k+1)t+k+2}{k}\\
>&0,
\end{align*}
which implies $|E(G_1)|<\binom{n-\left\lfloor\left(k+\frac{1}{2}\right)t\right\rfloor-1}{2}+t\left(\left\lfloor\left(k+\frac{1}{2}\right)t\right\rfloor+1\right)$.
Combining this with \eqref{eq:3.2}, we get
$|E(G)|\leq |E(G_1)|<\binom{n-\left\lfloor\left(k+\frac{1}{2}\right)t\right\rfloor-1}{2}+t\left(\left\lfloor\left(k+\frac{1}{2}\right)t\right\rfloor+1\right)$,
a contradiction. If $s\geq t+2$, then it follows from \eqref{eq:3.6} and $n\geq\left(k+\frac{3}{2}\right)s+1$ that
\begin{align}\label{eq:3.7}
&\binom{n-\left\lfloor\left(k+\frac{1}{2}\right)t\right\rfloor-1}{2}+t\left(\left\lfloor\left(k+\frac{1}{2}\right)t\right\rfloor+1\right)-|E(G_1)|\nonumber\\
\geq&\frac{1}{8}((2(2k+1)(s-t)+2)\cdot((2k+3)s+2)-(2k+1)(2k+5)s^{2}-(12k+14)s\nonumber\\
&+(2k+1)(2k+5)t^{2}+(8k+8)t-5)\nonumber\\
=&\frac{1}{8}((2k+1)^{2}s^{2}-(2(2k+1)(2k+3)t+4)s+(2k+1)(2k+5)t^{2}+4t-1)\nonumber\\
=&\frac{1}{8}\psi(s),
\end{align}
where $\psi(s)=(2k+1)^{2}s^{2}-(2(2k+1)(2k+3)t+4)s+(2k+1)(2k+5)t^{2}+4t-1$. Notice that
$$
\frac{2(2k+1)(2k+3)t+4}{2(2k+1)^{2}}<t+2\leq s
$$
by $1\leq t\leq k-1$. Then we obtain
\begin{align*}
\psi(s)\geq&\psi(t+2)\\
=&-8(2k+1)t+16k^{2}+16k-5\\
\geq&-8(2k+1)(k-5)+16k^{2}+16k-5\\
=&24k+3\\
>&0.
\end{align*}
Together with \eqref{eq:3.7}, we have
$$
\binom{n-\left\lfloor\left(k+\frac{1}{2}\right)t\right\rfloor-1}{2}+t\left(\left\lfloor\left(k+\frac{1}{2}\right)t\right\rfloor+1\right)-|E(G_1)|
\geq\frac{1}{8}\psi(s)>0,
$$
which yields $|E(G_1)|<\binom{n-\left\lfloor\left(k+\frac{1}{2}\right)t\right\rfloor-1}{2}+t\left(\left\lfloor\left(k+\frac{1}{2}\right)t\right\rfloor+1\right)$.
Together with \eqref{eq:3.2}, we get
$|E(G)|\leq |E(G_1)|<\binom{n-\left\lfloor\left(k+\frac{1}{2}\right)t\right\rfloor-1}{2}+t\left(\left\lfloor\left(k+\frac{1}{2}\right)t\right\rfloor+1\right)$,
a contradiction. This completes the proof of Theorem 1.3. \hfill $\Box$

\medskip

\noindent{\textbf{Remark 3.1.}} We are to show that the condition
$|E(G)|>\binom{n-\left\lfloor\left(k+\frac{1}{2}\right)t\right\rfloor-1}{2}+t\left(\left\lfloor\left(k+\frac{1}{2}\right)t\right\rfloor+1\right)$ in Theorem
1.3 cannot be replaced by
$|E(G)|\geq\binom{n-\left\lfloor\left(k+\frac{1}{2}\right)t\right\rfloor-1}{2}+t\left(\left\lfloor\left(k+\frac{1}{2}\right)t\right\rfloor+1\right)$.

We construct a graph
$G=K_t\vee(K_{n-\left\lfloor\left(k+\frac{3}{2}\right)t\right\rfloor-1}\cup\left(\left\lfloor\left(k+\frac{1}{2}\right)t\right\rfloor+1\right)K_1)$.
Then we have $|E(G)|=\binom{n-\left\lfloor\left(k+\frac{1}{2}\right)t\right\rfloor-1}{2}+t\left(\left\lfloor\left(k+\frac{1}{2}\right)t\right\rfloor+1\right)$.
Let $S=V(K_t)$. Then we obtain $|S|=t$. Thus, we deduce
\begin{align*}
i(G-S)=&\left\lfloor\left(k+\frac{1}{2}\right)t\right\rfloor+1\\
\geq&\left(k+\frac{1}{2}\right)t+\frac{1}{2}\\
=&\left(k+\frac{1}{2}\right)|S|+\frac{1}{2}\\
>&\left(k+\frac{1}{2}\right)|S|.
\end{align*}
According to Theorem 1.1, $G$ contains no $\{K_{1,1},K_{1,2},\ldots,K_{1,k},\mathcal{T}(2k+1)\}$-factor.

\section{The proof of Theorem 1.4}

\noindent{\it Proof of Theorem 1.4.} Suppose, to the contrary, that $G$ contains no $\{K_{1,1},K_{1,2},\ldots,K_{1,k},\mathcal{T}(2k+1)\}$-factor. Then
by Theorem 1.1, we conclude
$$
i(G-S)>\left(k+\frac{1}{2}\right)|S|
$$
for some subset $S$ of $V(G)$. In terms of the integrity of $i(G-S)$, we get
\begin{align}\label{eq:4.1}
i(G-S)\geq\left\lfloor\left(k+\frac{1}{2}\right)|S|\right\rfloor+1
\end{align}
for some subset $S$ of $V(G)$.

\noindent{\bf Claim 1.} $|S|\geq\delta$.

\noindent{\it Proof.} Assume that $|S|\leq\delta-1$. Since $G$ has minimum degree $\delta$, $G-S$ has minimum degree at least 1. Hence, we obtain
$i(G-S)=0$. Combining this with \eqref{eq:4.1}, we possess
$$
0=i(G-S)\geq\left\lfloor\left(k+\frac{1}{2}\right)|S|\right\rfloor+1\geq1,
$$
which is a contradiction. Claim 1 is proved. \hfill $\Box$

According to \eqref{eq:4.1} and Claim 1, we get
\begin{align}\label{eq:4.2}
i(G-S)\geq\left\lfloor\left(k+\frac{1}{2}\right)|S|\right\rfloor+1\geq\left\lfloor\left(k+\frac{1}{2}\right)\delta\right\rfloor+1.
\end{align}
From \eqref{eq:4.2}, we see that there exist at least $\left\lfloor\left(k+\frac{1}{2}\right)\delta\right\rfloor+1$ vertices
$x_1,x_2,\ldots,x_{\lfloor(k+\frac{1}{2})\delta\rfloor+1}$ in $I(G-S)$. Clearly, $d_{G-S}(x_i)=0$ for $1\leq i\leq\lfloor(k+\frac{1}{2})\delta\rfloor+1$.
Together with the condition of Theorem 1.4, we deduce
\begin{align}\label{eq:4.3}
\frac{2n}{2k+3}\leq&\max\{d_G(x_1),d_G(x_2),\ldots,d_G(x_{\lfloor(k+\frac{1}{2})\delta\rfloor+1})\}\nonumber\\
\leq&\max\{d_{G-S}(x_1),d_{G-S}(x_2),\ldots,d_{G-S}(x_{\lfloor(k+\frac{1}{2})\delta\rfloor+1})\}+|S|\nonumber\\
=&|S|.
\end{align}

It follows from \eqref{eq:4.1} and \eqref{eq:4.3} that
\begin{align*}
n\geq&|S|+i(G-S)\geq|S|+\left\lfloor\left(k+\frac{1}{2}\right)|S|\right\rfloor+1\\
=&\left\lfloor\left(k+\frac{3}{2}\right)|S|\right\rfloor+1\\
\geq&\left\lfloor\left(k+\frac{3}{2}\right)\cdot\frac{2n}{2k+3}\right\rfloor+1\\
=&n+1,
\end{align*}
which is a contradiction. This completes the proof of Theorem 1.4. \hfill $\Box$

\medskip

\noindent{\textbf{Remark 4.1.}} In what follows, we show that
$$
\max\{d_G(x_1),d_G(x_2),\ldots,d_G(x_{\lfloor(k+\frac{1}{2})\delta\rfloor+1})\}\geq\frac{2n}{2k+3}
$$
in Theorem 1.4 cannot be replaced by
$$
\max\{d_G(x_1),d_G(x_2),\ldots,d_G(x_{\lfloor(k+\frac{1}{2})\delta\rfloor+1})\}\geq\frac{2n-1}{2k+3}.
$$

We construct a graph $G=K_{\delta}\vee(\lfloor(k+\frac{1}{2})\delta\rfloor+1)K_1$ of order $n$, where $\delta\geq1$ is an odd integer. Obviously,
$\delta(G)=\delta$, $n=\lfloor(k+\frac{3}{2})\delta\rfloor+1=(k+\frac{3}{2})\delta+\frac{1}{2}$, and
$$
\max\{d_G(x_1),d_G(x_2),\ldots,d_G(x_{\lfloor(k+\frac{1}{2})\delta\rfloor+1})\}=\delta=\frac{2n-1}{2k+3}
$$
for any independent subset $\{x_1,x_2,\ldots,x_{\lfloor(k+\frac{1}{2})\delta\rfloor+1}\}$ of $G$. Set $S=V(K_{\delta})$. Then $|S|=\delta$ and
\begin{align*}
i(G-S)=&\left\lfloor\left(k+\frac{1}{2}\right)\delta\right\rfloor+1\\
=&\left(k+\frac{1}{2}\right)\delta+\frac{1}{2}\\
=&\left(k+\frac{1}{2}\right)|S|+\frac{1}{2}\\
>&\left(k+\frac{1}{2}\right)|S|.
\end{align*}
In terms of Theorem 1.1, $G$ contains no $\{K_{1,1},K_{1,2},\ldots,K_{1,k},\mathcal{T}(2k+1)\}$-factor.

\section{The proof of Theorem 1.5}

\noindent{\it Proof of Theorem 1.5.} Suppose, to the contrary, that $G$ has no $\{K_{1,1},K_{1,2},\ldots,K_{1,k},\mathcal{T}(2k+1)\}$-factor. Then
by Theorem 1.1, we obtain
\begin{align}\label{eq:5.1}
i(G-S)>\left(k+\frac{1}{2}\right)|S|
\end{align}
for some subset $S$ of $V(G)$. From \eqref{eq:5.1}, there exists at least one vertex $x$ in $I(G-S)$, and so $d_{G-S}(x)=0$. Thus, we deduce
$\delta(G)\leq d_G(x)\leq d_{G-S}(x)+|S|=|S|$. Together with \eqref{eq:5.1}, $\delta(G)\geq\frac{2\alpha(G)}{2k+1}$
and $\alpha(G)\geq i(G-S)$, we conclude
$$
|S|\geq\delta(G)\geq\frac{2\alpha(G)}{2k+1}\geq\frac{2i(G-S)}{2k+1}>\frac{2\left(k+\frac{1}{2}\right)|S|}{2k+1}=|S|,
$$
which is a contradiction. This completes the proof of Theorem 1.5. \hfill $\Box$

\medskip

\noindent{\textbf{Remark 5.1.}} We are to show that the condition $\delta(G)\geq\frac{2\alpha(G)}{2k+1}=\frac{\alpha(G)}{k+\frac{1}{2}}$ in Theorem
1.5 cannot be replaced by $\delta(G)\geq\frac{\alpha(G)-1}{k+\frac{1}{2}}$.

We construct a graph $G=K_{2+2t}\vee((1+t)(2k+1)+1)K_1$, where $t\geq0$
be an integer. Then we obtain $\delta(G)=2+2t$ and $\alpha(G)=(1+t)(2k+1)+1$. Thus, we deduce $\delta(G)=\frac{\alpha(G)-1}{k+\frac{1}{2}}$. Let
$S=V(K_{2+2t})$. Then we have $|S|=2+2t$. Thus, we conclude
\begin{align*}
i(G-S)=&(1+t)(2k+1)+1\\
=&\left(k+\frac{1}{2}\right)(2+2t)+1\\
=&\left(k+\frac{1}{2}\right)|S|+1\\
>&\left(k+\frac{1}{2}\right)|S|.
\end{align*}
By virtue of Theorem 1.1, $G$ contains no $\{K_{1,1},K_{1,2},\ldots,K_{1,k},\mathcal{T}(2k+1)\}$-factor.

\section*{Data availability statement}

My manuscript has no associated data.

\section*{Declaration of competing interest}

The authors declare that they have no known competing financial interests or personal relationships that could have appeared to influence the work
reported in this paper.

\section*{Acknowledgments}

This work was supported by the Natural Science Foundation of Jiangsu Province (Grant No. BK20241949). Project ZR2023MA078 supported by Shandong
Provincial Natural Science Foundation.

\end{document}